\documentclass[12pt]{article}
\usepackage{amsthm}
\usepackage{amsmath}
\usepackage{amsfonts}
\usepackage{epsfig}
\usepackage{subfigure}
\usepackage{a4wide}
\usepackage{amssymb}
\usepackage{amscd}
\usepackage{latexsym}
\newtheorem{theorem}{Theorem}
\newtheorem{lemma}[theorem]{Lemma}
\newtheorem{proposition}[theorem]{Proposition}
\newtheorem{corollary}[theorem]{Corollary}
\newtheorem{observation}[theorem]{Observation}
\newtheorem{conjecture}{Conjecture}

\newcommand{\qite}[1]{\noindent\leavevmode\hangindent1.5\parindent%
        \noindent\hbox to1.5\parindent{#1\hss}\ignorespaces}

\def\RR{{\mathbb R}}
\date{}
\begin{document}
\title{An improved linear bound on the number of perfect matchings in
cubic graphs\thanks{This research was partially supported by the
Czech-Slovenian bilateral project MEB 090805 and BI-CZ/08-09-005.}}
\author{Louis Esperet\thanks{Department of Applied Mathematics and
Institute for Theoretical Computer Science, Faculty of Mathematics and
Physics, Charles University, Malostransk{\'e} n{\'a}m{\v e}st{\'\i}
25, 118 00 Prague, Czech Republic.  E-mail: {\tt
esperet@kam.mff.cuni.cz}. Partially supported by the European
project \textsc{ist fet Aeolus}.}\and Daniel
Kr{\'a}l'\thanks{Institute for Theoretical Computer Science, Faculty
of Mathematics and Physics, Charles University, Malostransk\'e
n\'am\v{e}st\'{\i}~25, 118~00 Prague, Czech Republic.  E-mail: {\tt
kral@kam.mff.cuni.cz}.  The Institute for Theoretical Computer Science
is supported by Ministry of Education of the Czech Republic as project
1M0545.}\and Petr {\v Skoda}\thanks{Department of Applied Mathematics,
Faculty of Mathematics and Physics, Charles University, Malostransk\'e
n\'am\v{e}st\'{\i}~25, 118~00 Prague, Czech Republic.  E-mail: {\tt
peskoj@kam.mff.cuni.cz}.}\and Riste {\v S}krekovski\thanks{Department
of Mathematics, University of Ljubljana, Jadranska~19, 1111 Ljubljana,
Slovenia. Partially supported by ARRS Research Program
P1-0297.}}

\maketitle

\begin{abstract}
We show that every cubic bridgeless graph with $n$ vertices has at
least $3n/4-10$ perfect matchings. This is the first bound that
differs by more than a constant from the maximal dimension of the
perfect matching polytope.

%% This is the first bound for
%% cyclically $4$-edge-connected cubic graphs that differs by more than a
%% constant from the bound given by the dimension of the perfect matching
%% polytope.
\end{abstract}

\section{Introduction}
\label{sec-intro}

We study the number of perfect matchings in cubic bridgeless graphs,
in which parallel edges are allowed.  By a classical theorem of
Petersen~\cite{bib-petersen1891}, every such graph has a perfect
matching. In fact, every edge of a cubic bridgeless graph is contained
in a perfect matching, and thus every $n$-vertex cubic bridgless graph
has at least three perfect matchings.  Lov{\'a}sz and
Plummer~\cite[Conjecture 8.1.8]{bib-lovasz86+} conjectured that the
number of perfect matchings in cubic bridgeless graphs should grow
exponentially with $n$:

\begin{conjecture}[Lov{\'a}sz and Plummer, 1970s]
\label{conj}
Every cubic bridgeless graph with $n$ vertices
has at least $2^{\Omega(n)}$ perfect matchings.
\end{conjecture}

\noindent Conjecture~\ref{conj} has been verified for several special
classes of graphs, one of them being bipartite graphs.  The first
non-trivial lower bound on the number of perfect matchings in cubic
bridgeless bipartite graphs was obtained in 1969 by
Sinkhorn~\cite{bib-sinkhorn69} who proved a bound of $\frac{n}{2}$,
thereby establishing a conjecture of Hall. The same year,
Minc~\cite{bib-minc69} increased this lower bound by $2$.  Then, a
bound of $\frac{3n}{2}-3$ was proven by Hartfiel and
Crosby~\cite{bib-hartfiel71+}.  The first exponential bound,
$6\cdot\left(\frac{4}{3}\right)^{n/2-3}$, was obtained in 1979 by
Voorhoeve~\cite{bib-voorhoeve79}. This was generalized to all regular
bipartite graphs in 1998 by Schrijver~\cite{bib-schrijver98} who
thereby proved a conjecture of himself and
Valiant~\cite{bib-schrijver80+}.

Recently, an important step towards a proof of Conjecture~\ref{conj}
was achieved by Chudnovsky and Seymour~\cite{bib-chudnovsky08+}
who proved the conjecture for planar graphs.

\begin{theorem}[Chudnovsky and Seymour, 2008]
\label{thm-planar}
Every cubic bridgeless planar graph with $n$ vertices
has at least $2^{n/655978752}$ perfect matchings.
\end{theorem}

Until recently, the only known lower bound on the number of perfect
matchings of a general cubic bridgeless graph was an estimate given by
the dimension of the perfect matching polytope.  Edmonds, Lov\'asz,
and Pulleyblank~\cite{bib-edmonds82+}, inspired by
Naddef~\cite{bib-naddef82}, proved that the dimension of the perfect
matching polytope of a cubic bridgeless $n$-vertex graph is at least
$n/4+1$ which implies:

\begin{theorem}[Edmonds, Lov\'asz,
and Pulleyblank, 1982]
\label{thm-old}
Every cubic bridgeless graph with $n$ vertices has at least $n/4+2$
perfect matchings.
\end{theorem}

\begin{figure}
\begin{center}
\epsfbox{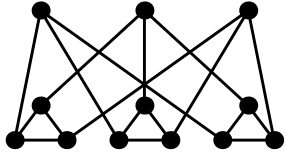}
\end{center}
\caption{The only $n$-vertex cubic bridgeless graph with $n/2$ perfect
matchings.}
\label{fig-exceptional}
\end{figure}

\noindent An argument based on the dimension of the perfect matching
polytope cannot yield a bound exceeding $n/2+2$, since the dimension
of the perfect matching polytope is always between $n/4+1$ and $n/2+1$
(the upper bound is achieved by cubic bipartite graphs).
In~\cite{bib-previous}, the authors presented an argument based on the
brick and brace decomposition of matching covered graphs, showing that
every $n$-vertex cubic bridgeless graph $G$ has at least $n/2$ perfect
matchings.  They also characterized those graphs $G$ with exactly
$n/2$ or $n/2+1$ perfect matchings. Their argument is inductive and
uses the characterization of so-called extremal cubic bricks by de
Carvalho \emph{et al.}~\cite{bib-carvalho05+}. Let us state the result
of~\cite{bib-previous} precisely:

\begin{theorem}
\label{thm-previous}
Every cubic bridgeless graph $G$ of order $n$ contains at least
$n/2+1$ perfect matchings unless $G$ is the graph obtained from
$K_{3,3}$ by replacing all three vertices of one of the two color
classes with triangles (see Figure~\ref{fig-exceptional}). This
exceptional graph contains $n/2$ perfect matchings.  Moreover, there
are only $17$ non-isomorphic cubic bridgeless graphs with at most
$n/2+1$ perfect matchings.
\end{theorem}

In this paper, we show that every $n$-vertex cubic bridgeless graph
has at least $3n/4-10$ perfect matchings.  We think that the main
significance of the bound lies in the fact that it is the first result
asserting that the number of vertices of the perfect matching polytope
of a cyclically $4$-edge-connected cubic graph exceeds its dimension
by more than a constant.

One of our tools, similarly as in~\cite{bib-previous}, is the
machinery of brick and brace decompositions of matching covered
graphs, which we introduce in the next section. However, unlike
in~\cite{bib-previous}, we have to show that the number of perfect
matchings of cyclically $4$-edge-connected graphs exceeds the
dimension of the perfect matching polytope by a linear factor.  This
is done in Sections~\ref{sect5} and \ref{sect4}. In
Sections~\ref{sect3} and \ref{sect2}, the bound is then extended to
$3$-edge-connected and eventually to all cubic bridgeless graphs.

\section{Notation}
\label{sect-defs}

In this section, we introduce notation used throughout the paper. If
$G$ is a graph, $V(G)$ denotes the vertex set of $G$ and $E(G)$
denotes its edge set. $\RR^{E(G)}$ is an $|E(G)|$-dimensional vector
space with coordinates corresponding to the edges of $G$. If
$A\subseteq V(G)$, $G[A]$ stands for the subgraph of $G$ induced by
the vertices of $A$.

A graph $G$ is {\em $k$-vertex-connected} if $G$ has at least $k+1$
vertices, and remains connected after removing any set of at most $k-1$
vertices. If $\{A,B\}$ is a partition of $V(G)$, the set $E(A,B)$ of
edges with one end in $A$ and the other in $B$ is called an {\em
edge-cut} or a {\it $k$-edge-cut} of $G$, where $k$ is the size of
$E(A,B)$. A graph is {\em $k$-edge-connected} if it has no edge-cuts
of size less than $k$. Graphs that are $2$-edge-connected are also
called {\em bridgeless}.  Finally, an edge-cut $E(A,B)$ is {\em
cyclic} if the subgraphs induced by $A$ and $B$ both contain a
cycle. A graph $G$ is {\em cyclically $k$-edge-connected} if $G$ has
no cyclic edge-cuts of size less than $k$. The following is a useful
observation that we implicitly use in our further considerations:

\begin{observation} If $G$ is a graph with minimum degree
three, in particular $G$ can be a cubic graph, then a $k$-edge-cut
$E(A,B)$ such that $|A|\ge k-1$ and $|B|\ge k-1$ must be cyclic.
\end{observation}

We say that a graph $G$ is {\it $X$-near cubic} for a multiset $X$ of
positive integers, if the multiset of degrees of $G$ not equal to
three is $X$. For example, the graph obtained from a cubic graph by
removing an edge is $\{2,2\}$-near cubic.

If $v$ is a vertex of $G$, then $G\setminus v$ is the graph obtained
by removing the vertex $v$ together with all its incident edges. If
$e$ is an edge of $G$, $G-e$ is the graph obtained from $G$ by
removing the edge $e$ and keeping its end vertices. We also use this
notation with $e$ replaced by a set of edges and $v$ replaced by a set
of vertices.  If $H$ is a connected subgraph of $G$, $G/H$ is the
graph obtained by contracting all the vertices of $H$ to a single
vertex, removing arising loops and preserving all parallel edges. An
{\em odd minor} of $G$ is a graph obtained by contracting connected
subgraphs of $G$, each having an odd number of vertices. Observe that
if all the degrees of $G$ are odd, then all the degrees of an odd
minor of $G$ are also odd.

A perfect matching of $G$ is a spanning subgraph with all vertices of
degree one.  A theorem of Tutte (1947) asserts that $G$ has a perfect
matching if and only if the number of components of $G\setminus S$
with an odd number of vertices (also called \emph{odd} components)
is at most $|S|$ for every $S\subseteq V(G)$. One of the consequences
of Tutte's theorem is that for every edge $e$ of a cubic bridgeless
graph, there is a perfect matching containing $e$ and for every two
edges $e$ and $f$, there is a perfect matching avoiding both $e$ and
$f$.

\subsection{Brick and brace decomposition of graphs}

The brick and brace decomposition plays a crucial role in the study of
the structure of perfect matchings in graphs.  A graph $G$ is said to
be {\em matching covered} if every edge is contained in a perfect
matching of $G$, and it is {\em matching double-covered} if every
edge is contained in at least two perfect matchings of $G$. A theorem
of Kotzig (see \cite[Section~8.6]{bib-lovasz86+}) asserts that if a
graph $G$ has a unique perfect matching, then $G$ has a bridge.  An
immediate consequence of this theorem is the following proposition:

\begin{proposition}
\label{prop-double}
Every cyclically $4$-edge-connected cubic graph different from $K_4$
is matching double-covered.
\end{proposition}

An edge-cut $E(A,B)$ is {\em tight} if every perfect matching contains
precisely one edge of $E(A,B)$. If $G$ is a connected matching covered
graph with a tight edge-cut $E(A,B)$, then $G[A]$ and $G[B]$ are also
connected. Moreover, every perfect matching of $G$ corresponds to a
pair of perfect matchings in the graphs $G/A$ and $G/B$. Hence, both
$G/A$ and $G/B$ are also matching covered. We say that we have
decomposed $G$ into $G/A$ and $G/B$. If any of these graphs still have
a tight edge-cut, we can keep decomposing it until no graph in the
decomposition has a tight edge-cut. Matching covered graphs without
tight edge-cuts are called \emph{braces} if they are bipartite and
\emph{bricks} otherwise, and the decomposition of a graph $G$ obtained
this way is known as the {\em brick and brace decomposition} of $G$.

Lov{\'a}sz~\cite{bib-lovasz87} showed that the collection of graphs
obtained from $G$ in any brick and brace decomposition is unique up to
the multiplicity of edges. This allows us to speak of \emph{the} brick
and brace decomposition of $G$, as well as \emph{the} number of bricks
and \emph{the} number of braces in the decomposition of $G$.

A graph is said to be {\em bicritical} if $G\setminus\{u,v\}$ has a
perfect matching for any two vertices $u$ and $v$. Edmonds \emph{et
al.}~\cite{bib-edmonds82+} gave the following characterization of
bricks:

\begin{theorem}[Edmonds \emph{et al.}, 1982]
\label{th-brick}
A graph $G$ is a brick if and only if it is $3$-vertex-connected and
bicritical.
\end{theorem}

It can also be proven that a brace is a bipartite graph such that for
any two vertices $u$ and $u'$ from the same color class and any two
vertices $v$ and $v'$ from the other color class, the graph
$G\setminus\{u,u',v,v'\}$ has a perfect matching,
see~\cite{bib-lovasz86+}.

We finish this subsection with an observation that the brick and brace
decomposition of a bipartite graph contains braces only; we include
the proof of this fact as a demonstration of the just introduced
notation.

\begin{proposition}
\label{prop-bip}
If $H$ is a bipartite matching covered graph, then its brick and
brace decomposition consists of braces only.
\end{proposition}

\begin{proof}
We proceed by induction on the size of $H$.  Let $U$ and $V$ be the
two color classes of $H$. If $H$ has no tight edge-cut, then $H$ is a
brace and there is nothing to prove. Otherwise, let $E(A,B)$ be a
tight edge-cut of $H$. Let $e$ be an edge of $E(A,B)$.  By symmetry,
we can assume that $e$ is incident with a vertex of $A\cap U$. Since
$H$ contains a perfect matching such that $e$ is the only edge of
$E(A,B)$ in the matching, $|A\cap U|=|A\cap V|+1$ and $|B\cap
V|=|B\cap U|+1$.  Hence, any matching containing a single edge of the
cut $E(A,B)$, say $f$, must satisfy that $f$ is incident with a vertex
of $A\cap U$. Since $E(A,B)$ is a tight edge-cut, all the edges of
$E(A,B)$ join vertices of $A\cap U$ and $B\cap V$, and so both graphs
$G/A$ and $G/B$ are bipartite. The claim follows by applying the
induction to $G/A$ and $G/B$.
\end{proof}

\subsection{Perfect matching polytope}

Some of our arguments also involve the perfect matching polytopes
of graphs.  The {\em perfect matching polytope} of a graph $G$ is the
convex hull of characteristic vectors of perfect matchings of $G$. The
sufficient and necessary conditions for a vector $w\in\RR^{E(G)}$ to
lie in the perfect matching polytope are known~\cite{bib-edmonds65}:

%% \begin{theorem}[Edmonds 1965]
%% \label{thm-polytope-bip}
%% If $G$ is a bipartite graph, then a vector $w\in\RR^{E(G)}$ lies in
%% the perfect matching polytope of $G$ if and only if $w$ is
%% non-negative and for every vertex $v$ of $G$ the sum of the entries of
%% $w$ corresponding to the edges incident with $v$ is equal to one.
%% \end{theorem}

\begin{theorem}[Edmonds 1965]
\label{thm-polytope}
If $G$ is a graph, then a vector $w\in\RR^{E(G)}$ lies in the perfect
matching polytope of $G$ if and only if the following holds:

\smallskip
  \qite{(i)}$w$ is non-negative,

  \smallskip
  \qite{(ii)}for every
vertex $v$ of $G$ the sum of the entries of $w$ corresponding to the
edges incident with $v$ is equal to one, and

  \smallskip \qite{(iii)}for every set $S \subseteq V(G)$, $|S|$ odd,
  the sum of the entries corresponding to edges having exactly one
  vertex in $S$ is at least one.
\end{theorem}

\noindent It is also well-known that conditions (i) and (ii) are necessary and
sufficient for a vector to lie in the perfect matching polytope of a
bipartite graph $G$.

The dimension of the perfect matching polytope of a matching covered
graph $G$ can be computed from the brick and brace decomposition of
$G$: Edmonds, Lov{\'a}sz, and Pulleyblank~\cite{bib-edmonds82+}, using
some ideas from Naddef~\cite{bib-naddef82}, showed that it is equal to
$|E(G)|-|V(G)|+1-b(G)$ where $b(G)$ denotes the number of bricks in
the decomposition.

Let $w$ be a vector lying in the perfect matching polytope of $G$ and
$E(A,B)$ be an edge-cut of $G$.  If the sum of the entries of $w$
corresponding to edges of $E(A,B)$ is not equal to one, then at least one of
the matchings whose characteristic vectors convexly combine to $w$
does not contain exactly one edge of the cut. Hence, $E(A,B)$ cannot
be tight. Conversely, if an edge-cut is tight, the entries
corresponding to the edges of the cut of every vector lying in the
perfect matching polytope sum to one. Let us formulate this
observation as a propostion.

\begin{proposition}
\label{prop-tight}
An edge-cut of $G$ is tight if and only if the sum of the entries
corresponding to the edges of the cut is equal to one for every vector
lying in the perfect matching polytope of $G$.
\end{proposition}

If $G$ is a cubic bridgeless graph, it is easy to infer from
Theorem~\ref{thm-polytope} that the vector with all entries equal to
$1/3$ lies in the perfect matching polytope of $G$.  Hence, every
tight cut of a cubic bridgeless graph must have size three by
Proposition~\ref{prop-tight}.  In particular, the brick and brace
decomposition of a cubic bridgeless graph only contains cubic
(bridgeless) graphs.

\section{Cyclically $5$-edge-connected graphs}
\label{sect5}

Our aim in this section is to show that if $G$ is a cyclically
$5$-edge-connected cubic graph, and $e$ is an edge of $G$, then $G-e$
has few bricks in its brick and brace decomposition, or there exists
an edge $f$ so that $G-\{e,f\}$ is bipartite and matching
covered. This will imply that $G$ has at least $3|V(G)|/4-3/2$ perfect
matchings.

\begin{lemma}
\label{lm-e-5-conn}
Let $G$ be a cyclically $5$-edge-connected cubic graph, and let
$E(U,U')$ be a 5-edge-cut of $G$. If $G/U$ is matching covered, then
it is cyclically $5$-edge-connected and $3$-vertex-connected.

%% Let $G$ be a cyclically $5$-edge-connected cubic graph, and let
%% $E(U,U')$ be a 5-edge-cut such that $G[U]$ is connected. If $G/U$ is
%% matching covered, then it is cyclically $5$-edge-connected and
%% $3$-vertex-connected.
\end{lemma}

\begin{proof}
Since $G$ is cyclically $5$-edge-connected, $G[U]$ is connected, and
so $H=G/U$ is well-defined. Observe that any cyclic edge-cut of $H$
corresponds to a cyclic edge-cut of $G$. Hence, $H$ is cyclically
$5$-edge-connected. Moreover, it is a $\{5\}$-near cubic graph, and
since the minimum degree of $H$ is three, any edge-cut of size at most
two is cyclic. This implies that $H$ is 3-edge-connected. Also note
that $H$ is 2-vertex-connected, otherwise it would contain an edge-cut
of size at most two since the maximum degree of $H$ is five.

We now show that $H$ is 3-vertex-connected, which will establish the
lemma. For the sake of contradiction, assume that $H$ has a vertex-cut
of size two formed by vertices $x$ and $y$, and let $A$ and $B$ be the
components of $H\setminus\{x,y\}$. If both $x$ and $y$ have degree
three, one easily infer a 2-edge-cut.
%% the number of edges between $A$ (resp. $B$) and $\{x,y\}$ is
%% equal to three since $H$ is 3-edge-connected. If $A$ (resp. $B$)
%% contains the vertex of degree five, then the edges between $A$
%% (resp. $B$) and $\{x,y\}$ form a cyclic edge-cut of $G$ which is
%% impossible.
Hence, we may assume that $x$ has degree five and $y$ has
degree three. By the 3-edge-connectivity of $H$,
the graph $H \setminus \{x,y\}$ cannot have more than two
components.

%Note that the number of vertices of $H$ is even. So, if $H$ has order
%four, then the underlying simple graph of $H$ must be $C_4$ as $H$ is
%matching covered. A simple check implies that $H$ cannot be
%$\{5\}$-near cubic in this case. Hence, $H$ has at least six vertices. 

A simple check shows that the only $\{5\}$-near cubic graph of order
at most four is the graph obtained from $K_4$ by removing an edge, say
$uv$, and doubling the edges $uw$ and $vw$, where $w$ is one of the
two vertices distinct from $u$ and $v$. However, this graph is not
matching covered. Since the number of vertices of $H$ is even, we can
assume that $H$ has at least six vertices.

If $x$ and $y$ are joined by an edge, then the number of edges between
$A$ or $B$ and $\{x,y\}$ must be three. At least one these two
edge-cuts is however cyclic; otherwise, both $A$ and $B$ have order
one and the order of $H$ is four. Hence, the number of edges leaving
$\{x,y\}$ is eight and $x$ and $y$ are non-adjacent.

Neither $x$ nor $y$ is incident with a bigon (an edge with
multiplicity two); otherwise the edges leaving the bigon form a cyclic
edge-cut of $H$ of size at most four. Since the number of edges
between $A$ or $B$ and $\{x,y\}$ must be at least three and neither
$x$ nor $y$ is incident with a bigon, it follows that both $A$ and $B$
contain at least two vertices. Hence, the number of edges between $A$
or $B$ and $\{x,y\}$ must be at least four since otherwise these edges
would form a cyclic edge-cut of size three in $H$. Consequently, there
are exactly four edges between $A$ or $B$ and $\{x,y\}$, and the sets
$A$ and $B$ both contain exactly two vertices. Since $x$ has degree
five and is neither adjacent to $y$ nor incident to a bigon, this is
impossible.
\end{proof}

We now prove that under the same assumptions as in the previous lemma,
the brick and brace decomposition of $G/U$ contains exactly one brick.

\begin{lemma}
\label{lm-e-5}
Let $G$ be a cyclically $5$-edge-connected cubic graph, and let
$E(U,U')$ be a 5-edge-cut of $G$. If $G/U$ is matching covered, then
$b(G/U)=1$.

%% Let $G$ be a cyclically $5$-edge-connected cubic graph, and let
%% $E(U,U')$ be a 5-edge-cut such that $G[U]$ is connected. If $G/U$ is
%% matching covered, then $b(G/U)=1$.
\end{lemma}

\begin{proof}
The proof proceeds by induction on the order of $H=G/U$ ($G$ is
fixed). Since $H$ is a $\{5\}$-near cubic graph, $H$ is not
bipartite. By Lemma~\ref{lm-e-5-conn}, $H$ is cyclically
$5$-edge-connected and $3$-vertex-connected. By
Theorem~\ref{th-brick}, $H$ is a either a brick (in which case
$b(H)=1$) or is not bicritical. So we can focus on the latter case.

Let $x$ and $y$ be the vertices of $H$ such that $H\setminus\{x,y\}$
has no perfect matching.  According to Tutte's Theorem, there exists a
set of vertices $S$ of $H\setminus\{x,y\}$ such that
$H\setminus(S\cup\{x,y\})$ has at least $|S|+1$ odd components. Let
$S'=S\cup\{x,y\}$.  Since the number of vertices of $H$ is even,
$H\setminus S'$ has at least $|S|+2=|S'|$ odd components. As $H$
is $\{5\}$-near cubic, the number of edges leaving $S'$ is at most
$3|S'|+2$. In what follows, we distinguish two cases regarding the
sizes of the components in $H\setminus S'$.

Suppose first that all the components of $H\setminus S'$ are single
vertices of degree three in $H$. Then the number of edges between $S'$
and $H\setminus S'$ is exactly $3|S'|$. In this case, the vertex of
degree five is in $S'$ and $S'$ contains two vertices joined by an
edge. Observe that $H$ has no matching containing this edge which
contradicts our assumption that $H$ is matching covered.

Suppose now that at least one of the components of $H\setminus S'$ is
not a single vertex whose degree is three in $H$, then the number of
edges leaving the odd components of $H\setminus S'$ is at least
$3|S'|+2$: there are at least five edges leaving every odd component
that is not a single vertex since $H$ is cyclically $5$-edge-connected
and there are five edges leaving a vertex of degree five in case this
vertex were one of the components of $H\setminus S'$.  We conclude
that the number of edges between $S'$ and $H\setminus S'$ is exactly
$3|S'|+2$ (and thus $S'$ is a stable set and contains the vertex of
degree five), and $H\setminus S'$ contains exactly $|S'|$ components,
$|S'|-1$ of them being isolated vertices and the remaining one having
odd size.

Let $B$ be the set of vertices of the only component of $H\setminus
S'$ that is not an isolated vertex and set $A=V(H)\setminus B$. As
$H\setminus S'$ contains exactly $|S'|$ components and $S'$ is a
stable set, the 5-edge-cut $E(A,B)$ is tight. In particular, $H/B$ is
a bipartite matching covered graph, so $b(H/B)=0$ by
Proposition~\ref{prop-bip}. Let $A'$ be the set of vertices of $G$
corresponding to $A$, i.e. $H/A=G/A'$. The graph $H[A]$
is connected and contains the vertex of degree five, so $H/A=G/A'$ is
a matching covered graph that satisfies the induction
hypothesis. Since the order of $G/A'=H/A$ is smaller than that of
$G/U=H$, the induction yields that $b(H/A)=1$. Consequently,
$b(H)=b(H/A)+b(H/B)=1+0=1$.
\end{proof}

Using the same approach as in Lemma~\ref{lm-e-5-conn}, we now study
the connectivity of a matching covered $\{4,4\}$-near cubic graph
obtained from $G-e$ by contracting some odd components.

\begin{figure}
\begin{center}
\epsfbox{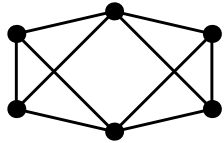}
\end{center}
\caption{The only possible odd minor of $G-e$ ($G$ being a cyclically
         $5$-edge-connected cubic graph) that is a matching covered
         $\{4,4\}$-near cubic graph and that is not
         $3$-vertex-connected.}
\label{fig-e-44}
\end{figure}

\begin{lemma}
\label{lm-e-44-conn}
Let $G$ be a cyclically $5$-edge-connected cubic graph, $e$ an edge of
$G$ and $H$ an odd minor of $G-e$. If $H$ is a $\{4,4\}$-near cubic
graph, then $H$ is $2$-vertex-connected. Moreover, if $H$ has a
$2$-vertex-cut and is matching covered, then $H$ is isomorphic to the
graph depicted in Figure~\ref{fig-e-44}.
\end{lemma}

\begin{proof}
Since $G$ is cyclically $5$-edge-connected, the graph $H$ is
cyclically $4$-edge-connected. We first show that we can focus on
graphs $H$ of order six or more.  The only $\{4,4\}$-near cubic graphs
of order at most four that are matching covered but not
$3$-vertex-connected have the cycle $C_4$ as an underlying simple
graph. In that case, $H$ must be either
\begin{itemize}
\item[$\bullet$] the graph obtained from $C_4$ by doubling three
distinct edges, or
\item[$\bullet$] the graph obtained from $C_4$ by tripling an edge and
doubling the opposite one.
\end{itemize}
Since both these graphs contain a cyclic edge-cut of size at most
three, the order of $H$ is at least six.

The argument that $H$ is $2$-vertex-connected is analogous to that in
the proof of Lemma~\ref{lm-e-5-conn}, so we leave the details to the
reader. Assume that $H$ contains a bigon. Since $H$ has order at least
six and is cyclically $4$-edge-connected, exactly four edges leave
this bigon. Observe that $e$ is not contained in the corresponding
edge-cut in $G$, since the ends of the bigon are the two vertices of
degree four of $H$. Hence, the four edges leaving the bigon correspond
to a cyclic 4-edge-cut of $G$, which is impossible. So we can assume
that $H$ is a simple graph.

%% In addition, $H$ is a simple graph since otherwise the edges leaving a
%% bigon form a cyclic edge-cut of $H$ of size at most four. Note that in
%% case that there is a bigon joining the two vertices of degree four
%% (which is the only case when there could be four edges leaving the
%% bigon), the edge $e$ does not cross this cyclic edge-cut since it
%% corresponds to an edge joining the two vertices of degree four of
%% $H$. Since $e$ corresponds to an edge joining these two vertices, the
%% four edges leaving the bigon form cyclic 4-edge-cut of $G$, which is a
%% contradiction.

%Finally, we focus on analyzing vertex-cuts of size two. Let $\{x,y\}$
%be a $2$-vertex-cut of $H$ and let $A$ and $B$ be the two components
%of $H\setminus\{x,y\}$. Since $H$ has no bigons, each of the sets $A$
%and $B$ contains at least two vertices. Hence, the number of edges
%between $A$ and $\{x,y\}$ is at least four (otherwise, they would form
%a cyclic edge-cut of $H$ of size at most three). Similarly, the number
%of edges between $B$ and $\{x,y\}$ is at least four. Consequently, $x$
%and $y$ are non-adjacent vertices of degree four, the number of edges
%between $A$ (resp. $B$) and $\{x,y\}$ is precisely four.  Since the
%edge $e$ corresponds to an edge joining $x$ and $y$, each of the cuts
%$E(A,B\cup\{x,y\})$ and $E(A\cup\{x,y\},B)$ has the same size in $H$
%and $G$. As $G$ is cyclically $5$-edge-connected, both $A$ and $B$
%must contain exactly two vertices. A simple analysis using the fact
%that $H$ is simple yields that $H$ is the graph depicted in
%Figure~\ref{fig-e-44}.

Finally, we focus on analyzing vertex-cuts of size two. Let $\{x,y\}$
be a $2$-vertex-cut of $H$ and let $C_1,\ldots,C_k$ be the $k\ge 2$
components of $H\setminus\{x,y\}$. Since $H$ has no bigons, each of
the sets $C_i$ contains at least two vertices. Hence, the number of
edges between $C_i$ and $\{x,y\}$ is at least four (otherwise, they
would form a cyclic edge-cut of $H$ of size at most
three). Consequently, $k=2$ and $x$ and $y$ are non-adjacent vertices
of degree four. This implies that the number of edges between each
$C_i$ ($i=1,2$) and $\{x,y\}$ is precisely four. Since the edge $e$
corresponds to an edge joining $x$ and $y$, each of the cuts
$E(C_1,C_2\cup\{x,y\})$ and $E(C_1\cup\{x,y\},C_2)$ has the same size
in $H$ and $G$. As $G$ is cyclically $5$-edge-connected, both $C_1$
and $C_2$ must contain exactly two vertices. We conclude that $H$ must
be the graph depicted in Figure~\ref{fig-e-44}.
\end{proof}

In the next lemma, we show that graphs satisfying the assumptions of
Lemma~\ref{lm-e-44-conn} have few bricks in their decomposition.

\begin{lemma}
\label{lm-e-44}
Let $G$ be a cyclically $5$-edge-connected cubic graph, $e$ an edge of
$G$ and $H$ an odd minor of $G-e$. If $H$ is a matching covered
$\{4,4\}$-near cubic graph, then $b(H) \le 2$.
\end{lemma}

\begin{proof}
The proof proceeds by induction on the order of the odd minor $H$ of
$G-e$ ($G$ and $e$ are fixed). If $H$ is bipartite, then $b(H)=0$ by
Proposition~\ref{prop-bip}.  If $H$ is not $3$-vertex-connected, then
by Lemma~\ref{lm-e-44-conn} it is isomorphic to the graph depicted in
Figure~\ref{fig-e-44} and its brick and brace decomposition consists
of two graphs isomorphic to $K_4$ with a single parallel edge.

Hence, we can assume that $H$ is a 3-vertex-connected non-bipartite
graph. If $H$ is a brick, then $b(H)=1$. By Theorem~\ref{th-brick}, we
may assume that $H$ is not bicritical. Let $x$ and $y$ be two vertices
of $H$ such that $H\setminus\{x,y\}$ has no perfect matching.  Let $S$
be a set of vertices of $H\setminus\{x,y\}$ such that
$H\setminus(S\cup\{x,y\})$ has at least $|S|+1$ odd components, and
let $S'=S\cup\{x,y\}$. Since the number of vertices of $H$ is even,
$H\setminus S'$ has at least $|S|+2=|S'|$ odd components. Based on the
degree distribution of $H$ and the fact that $G$ is cyclically
$5$-edge-connected, the number of edges leaving $S'$ is $3|S'|$,
$3|S'|+1$ or $3|S'|+2$ and $H\setminus S'$ contains precisely $|S'|$
components (which are all odd) and at most one of these components is
not an isolated vertex. Notice that if all the odd components of
$H\setminus S'$ were isolated vertices, then either $H$ would be
bipartite (which case has already been considered) or $S'$ would
contain both vertices of degree four. In the latter case, there would
be an edge joining two vertices of $S'$ but such an edge cannot be
contained in a perfect matching of $H$ contrary to our assumption that
$H$ is matching covered. We conclude that $H\setminus S'$ contains
precisely one non-trivial odd component $B$.

Let $A=V(H)\setminus B$. We consider three possibilities, regarding
whether the vertices of degree four belong to $S'$. If $S'$ only
contains vertices of degree three, then there are $3|S'|$ edges
leaving $S'$. In this case, the two vertices of degree four are in $B$
and $E(A,B)$ is a cyclic edge-cut of size three, which is
impossible. Depending whether $S'$ contains one or both vertices of
degree four of $H$, the number of edges between $A$ and $B$ is four or
five. Observe that in both cases, $H/B$ is bipartite, and hence the
edge-cut $E(A,B)$ is tight. By Proposition~\ref{prop-bip}, this also
implies that $b(H/B)=0$. Let $A'$ the set of vertices of $G$
corresponding to $A$, i.e. $H/A=G/A'$. Since $E(A,B)$ is tight, the
graph $H/A=(G-e)/A'$ is matching covered. If $S'$ contains a single
vertex of degree four, then $H/A$ is a $\{4,4\}$-near cubic graph. In
this case we apply induction on $H/A$. If $S'$ contains two vertices
of degree four, then $E(A,B)$ is a cyclic 5-edge-cut and we can apply
Lemma~\ref{lm-e-5} on $H/A$. In both cases, $b(H)=b(H/A)+b(H/B)=
b(H/A)\le 2$.
\end{proof}

Lemma~\ref{lm-e-44} has the following corollary:

\begin{lemma}
\label{lm-e}
Let $G$ be a cyclically $5$-edge-connected cubic graph and $e$ an edge
of $G$. If $G-e$ is matching covered, then $b(G-e)\le 2$.
\end{lemma}

\begin{proof}
Since $G$ is a cyclically $5$-edge-connected cubic graph and $G-e$ is
matching covered, we infer that $G$ is not isomorphic to $K_4$. This
implies that $G$ is triangle-free. Hence, the two vertices of degree
two of $G-e$, say $u$ and $u'$, have no common neighbor.

Let $A$ be comprised of the vertex $u$ and its two neighbors in $G-e$
and $B=V(G)\setminus A$. Similarly, let $A'$ be comprised of the
vertex $u'$ and its two neighbors in $G-e$ and $B'=V(G)\setminus
A'$. The cuts $E(A,B)$ and $E(A',B')$ are tight in $G-e$. Since the
sets $A$ and $A'$ are disjoint, after reducing the tight edge-cuts
$E(A,B)$ and $E(A',B')$ of $G-e$, we obtain two bipartite graphs of
order four and a $\{4,4\}$-near cubic graph. The statement follows
from Proposition~\ref{prop-bip} and Lemma~\ref{lm-e-44}.
\end{proof}

We now study the structure of a graph $G$ such that the graph $G-e$ is
not matching covered for some edge $e$.

\begin{lemma}
\label{lm-e-non}
Let $G$ be a cyclically $5$-edge-connected cubic graph and $e$ an edge
of $G$.  If $G-e$ is not matching covered, then $G$ contains an edge
$f$ such that $G-\{e,f\}$ is matching covered and bipartite.
\end{lemma}

\begin{proof}
Let $e=uu'$ and $H=G-e$, and assume that $H$ contains an edge $f=vv'$
that is not contained in any perfect matching of $H$. Hence,
$H\setminus \{v,v'\}$ contains a set $S$ of vertices such that the
number of odd components of $H\setminus S'$ where $S'=S\cup\{v,v'\}$
is at least $|S|+1$. Since the number of vertices of $H$ is even, the
number of the odd components is at least $|S|+2=|S'|$. Since $v$ and
$v'$ are both contained in $S'$, the number of edges leaving $S'$ is
at most $3|S'|-2$. Since $G$ is cyclically $5$-edge-connected, all the
components of $H\setminus S'$ are isolated vertices and neither $u$
nor $u'$ is contained in $S'$.  This implies that
$H'=G\setminus\{e,f\}$ is a $\{2,2,2,2\}$-near cubic bipartite graph.
Denote by $U$ and $V$ the two color classes of $H'$, in such way that
$\{u,u'\}\subseteq U$ and $\{v,v'\}\subseteq V$.

We now show that $H'$ is matching covered. Let $H''$ be a graph
obtained from $H'$ by adding a vertex $v_e$ (resp. $v_f$) and joining
it by two parallel edges to each of the end-vertices of $e$
(resp. $f$). We claim that $H''$ has no edge-cut of size at most three
separating $v_e$ and $v_f$. Assume the opposite and let $E(A,B)$ be
such an edge-cut. By symmetry, $v_e\in A$ and $v_f\in B$.

If $A$ contains both end-vertices of $e$ and $B$ contains both
end-vertices of $f$, then $E(A,B)$ corresponds to a non-trivial
edge-cut of size at most three of $G$ which violates our assumption
that $G$ is cyclically $5$-edge-connected. Hence, we can assume by
symmetry that $A$ contains $u$ but not $u'$. As the size of $E(A,B)$ is
at most three, both $v$ and $v'$ must be contained in $B$. Let us
estimate the size of the edge-cut of $G$ corresponding to $E(A,B)$: the
two edges between $v_e$ and $u'$ are not present anymore and but the
edge $e$ is now present. Hence, the size of the corresponding edge-cut
of $G$ is at most two. Since $G$ is cubic, this is also a cyclic
edge-cut of size at most two, which contradicts our assumption that
$G$ is cyclically $5$-edge-connected.

Since there is no edge-cut of size at most three separating $v_e$ and
$v_f$ in $H''$, there are four edge-disjoint paths connecting $v_e$
and $v_f$ by Menger's theorem. Consequently, $H'$ contains four
edge-disjoint paths $P_1$, $P_2$, $P_3$ and $P_4$ joining the vertices
$u$ and $u'$ to the vertices $v$ and $v'$. Direct the paths $P_i$ from
$u$ and $u'$ to $v$ and $v'$, and consider now the following vector
$w\in\RR^{E(H')}$:
$$
w_e=\left\{\begin{array}{cl}
1/2 & \mbox{if $e$ is directed from $U$ to $V$,}\\
1/6 & \mbox{if $e$ is directed from $V$ to $U$, and}\\
1/3 & \mbox{otherwise.}
\end{array}\right.
$$ Observe that $H'$ is bipartite and for every vertex $x$ of $H'$,
the sum of the entries of $w$ corresponding to the edges incident with
$x$ is equal to one. Hence, $w$ lies in the perfect matching polytope
of $H'$. Since all the entries of $w$ are non-zero, the graph $H'$ is
matching covered.
\end{proof}

We now apply Lemmas~\ref{lm-e} and~\ref{lm-e-non} to prove the main
result of this section.

\begin{theorem}
\label{thm-5conn2}
Let $G$ be a cyclically $5$-edge-connected cubic graph of order
$n$. For every edge $e$ of $G$, the graph $G-e$ has at least $n/2-1$
perfect matchings.
\end{theorem}

\begin{proof}
Let $e$ be an arbitrary edge of $G$. If $G-e$ is matching covered,
then $b(G-e) \le 2$ by Lemma~\ref{lm-e}. Hence, the dimension of the
perfect matching polytope of $G-e$ is at least
$(3n/2-1)-n+1-2=n/2-2$. Consequently, $G-e$ has at least $n/2-1$
perfect matchings.

If $G-e$ is not matching covered, then Lemma~\ref{lm-e-non} guarantees
the existence of an edge $f$ such that $G\setminus\{e,f\}$ is matching
covered and bipartite, in which case $b(G\setminus\{e,f\})=0$ by
Proposition~\ref{prop-bip}. Hence, the dimension of the perfect
matching polytope of $G\setminus\{e,f\}$ is at least
$(3n/2-2)-n+1=n/2-1$ and $G-e$ contains at least $n/2$ perfect
matchings.
\end{proof}

This theorem has the following easy consequence on the number of
perfect matchings of cyclically $5$-edge-connected cubic graphs.

\begin{corollary}
\label{thm-5conn}
Let $G$ be a cubic graph of order $n$. If $G$ is cyclically
$5$-edge-connected, then the number of perfect matchings of $G$ is at
least $3n/4-3/2$.
\end{corollary}

\begin{proof}
Let $e$, $e'$ and $e''$ be the edges incident with an arbitrary vertex
$v$. By Theorem~\ref{thm-5conn2}, each of the graphs $G-e$, $G-e'$ and
$G-e''$ has at least $n/2-1$ perfect matchings. Since a perfect
matching of $G$ is a perfect matching of exactly two of these three
graphs, $G$ has at least $3n/4-3/2$ perfect matchings.
\end{proof}

\section{Cyclically $4$-edge-connected graphs}
\label{sect4}

In this section, we prove that cyclically $4$-edge-connected cubic
graphs have at least $3n/4-9$ perfect matchings. Actually, we prove a
slightly stronger version of this result that will be used in the next
section.

\begin{theorem}
\label{thm-4conn}
Let $H$ be a cyclically $4$-edge-connected cubic graph that is not
cyclically $5$-edge-connected. If $G$ is a graph of order $n$ obtained
from $H$ by replacing some of its vertices with triangles (possibly,
$G=H$), then $G$ contains at least $3n/4-9$ perfect matchings.
\end{theorem}

\begin{proof}
Let $E(A',B')=\{e'_1,e'_2,e'_3,e'_4\}$ be a cyclic $4$-edge-cut of
$H$. Let $a'_i$ be the end-vertex of the edge $e'_i$ lying in
$A'$. Observe that all the vertices $a'_i$ are distinct, since
otherwise there would be a cyclic edge-cut of size at most three in
$H$. We claim that the graph $H[A']$ is connected and bridgeless: If
$H[A']$ were disconnected, then a proper subset of
$\{e'_1,e'_2,e'_3,e'_4\}$ would also be a cyclic edge-cut which is
impossible by our assumption that $H$ is cyclically
$4$-edge-connected. If $H[A']$ has a bridge $e'$, this bridge must
separate in $A'$ two of the vertices $a'_1,a'_2,a'_3,a'_4$ from the
other two (otherwise, $H$ would contain an edge-cut of size
two). Assume that the bridge $e'$ separates $\{a'_1,a'_2\}$ from
$\{a'_3,a'_4\}$. As $\{e',e'_1,e'_2\}$ is an edge-cut of $H$ of size
three, $a'_1$ and $a'_2$ must coincide (otherwise, this edge-cut is
cyclic).  Similarly, we infer that $a'_3=a'_4$. This implies that the
subgraph $H[A']$ is just an edge contrary to the fact that $E(A',B')$
is a cyclic edge-cut. Hence, $H[A']$ and $H[B']$ are 2-edge-connected.

Observe that $E(A',B')$ corresponds to a cyclic $4$-edge-cut
$E(A,B)=\{e_1,e_2,e_3,e_4\}$ of $G$. Let $a_i$ and $b_i$ be the
end-vertex of the edge $e_i$ lying in $A$ and $B$, respectively. Now,
let $m^A_X$, $X\subseteq\{1,2,3,4\}$, be the number of matchings of
$G[A]$ that cover all the vertices of $G[A]$ except the vertices
$a_i$, $i\in X$. We use $m^B_X$ in an analogous way. To simplify our
notation, we further write $m^A_{13}$ instead of $m^A_{\{1,3\}}$,
etc. Clearly, if $|X|$ is odd, then $m^A_X=m^B_X=0$. As the number of
matchings of $G$ is equal to $$\sum_{X\subseteq\{1,2,3,4\}} m^A_X\cdot
m^B_X\;\mbox{,}$$ we will estimate the summands to obtain the desired
bound. Consider a permutation $\{i,j,k,l\}$ of $\{1,2,3,4\}$ with
$i<j$, and define $G^A_{ij}$ as the graph obtained from $G[A]$ by
adding the edges $a_ia_j$ and $a_ka_l$. $G^A_{(ij)}$ denotes the graph
obtained from $G[A]$ by introducing two new adjacent vertices, joining
one of them to the vertices $a_i$ and $a_j$, and the other one to
$a_k$ and $a_l$. Observe that $G^A_{12}=G^A_{34}$ and
$G^A_{(12)}=G^A_{(34)}$.

Since $H[A']$ is 2-edge-connected, so is the graph $G[A]$. Hence, the
graphs $G^A_{ij}$ and $G^A_{(ij)}$ are cubic and
bridgeless. Consequently, they have a perfect matching containing any
prescribed edge and a perfect matching avoiding any two prescribed
edges. In particular, $G^A_{12}$ has a matching avoiding the edges
$a_1a_2$ and $a_3a_4$. Consequently, $G[A]$ has a perfect
matching. Since $G[A]$ is bridgeless, it has at least two perfect
matchings by Kotzig's theorem. We conclude that $m^A_{\varnothing} \ge
2$. Also by Theorem~\ref{thm-previous}, the graphs $G^A_{ij}$ have at
least $|A|/2$ perfect matchings and the graphs $G^A_{(ij)}$ have at
least $|A|/2+1$ perfect matchings.

If $m^A_{1234}=0$, then the fact that $G^A_{ij}$ has a perfect
matching containing the edge $a_ia_j$ implies that $m^A_{ij}\ge 1$ for
every $i,j$. On the other hand, if $m^A_{ij}=0$ for some $i,j$ and
$k\not\in \{i,j\}$, then the fact that $G^A_{(jk)}$ has a perfect
matching containing the added edge incident with $a_i$ implies that
$m^A_{ik}\ge 1$. We conclude that at least one of the following two
possibilities occurs:
\begin{description}
\item[Case A:] \emph{All the quantities $m^A_{ij}$ are non-zero and
$m^A_\varnothing\ge 2$.}
\item[Case B:] \emph{There exist $i$ and $j$ such that the quantities
              $m^A_{1234}$, $m^A_{ik}$ and $m^A_{jk}$ are non-zero for
              any $k\not\in\{i,j\}$, and $m^A_\varnothing \ge 2$.}
\end{description}
For every subset $X\subseteq\{1,2,3,4\}$ such that $m^A_X\ge 1$, fix a
matching $M^A_X$ avoiding the vertices $a_i$, $i\in X$. In addition,
fix a second matching $M^{A*}_\varnothing\ne M^{A}_\varnothing$ covering
all the four vertices $a_i$, $i\in\{1,2,3,4\}$ (such a matching
exists as $m^A_\varnothing \ge 2$). The fixed matchings of $G[A]$ are
referred to as {\em canonical} matchings of $G[A]$ and the other
matchings of $G[A]$ are {\em non-canonical}. Consider also
the analogous definitions for the matchings of $G[B]$.

Assume first that Case A applies. Consider a non-canonical matching of
$G[B]$ that avoids vertices $b_i$ and $b_j$ for some
$i,j\in\{1,2,3,4\}$. This matching can be completed by adding the
canonical matching $M^A_{ij}$ and the edges $a_ib_i$ and $a_jb_j$ to a
perfect matching of $G$. Similarly, a non-canonical matching of $G[B]$
covering all the four vertices can be completed by one of the two
canonical matchings $M^{A}_\varnothing$ and $M^{A*}_\varnothing$ of
$G[A]$. We conclude that the number of perfect matchings of $G$ that
are canonical when restricted to $G[A]$ and non-canonical when
restricted to $G[B]$ is at least
\begin{equation}
\overline{m}^B_{12} +\overline{m}^B_{13} +\overline{m}^B_{14}
+\overline{m}^B_{23} +\overline{m}^B_{24} +\overline{m}^B_{34}
+2\overline{m}^B_{\varnothing},
\label{eq-4conn-1a}
\end{equation}
where $\overline{m}^B_X$ denotes the number of non-canonical matchings
of $G[B]$ avoiding $\{b_i, i\in X\}$. On the other hand, if
$\{i,j,k,l\}$ is a permutation of $\{1,2,3,4\}$, the number of perfect
matchings of $G^B_{(ij)}$ is equal to
\begin{equation}
m^B_{ik}+m^B_{il}+m^B_{jk}+m^B_{jl}+m^B_{\varnothing}.
\label{eq-4conn-2a}
\end{equation}
Every graph $G^B_{(ij)}$ has order $|B|+2$, so the number of perfect
matchings of $G^B_{(ij)}$ is at least $|B|/2+1$ by
Theorem~\ref{thm-previous} (and thus the number of non-canonical
matchings of $G[B]$ is at least $|B|/2-5$). Summing
(\ref{eq-4conn-2a}) for $(i,j)\in
\{(1,2),(1,3),(1,4)\}$ yields the following estimate:
\begin{equation}
2\overline{m}^B_{12} +2\overline{m}^B_{13} +2\overline{m}^B_{14}
+2\overline{m}^B_{23} +2\overline{m}^B_{24} +2\overline{m}^B_{34}
+3\overline{m}^B_{\varnothing}\ge 3|B|/2-15.
\label{eq-4conn-3a}
\end{equation}
Comparing (\ref{eq-4conn-1a}) and (\ref{eq-4conn-3a}), we see that the
number of perfect matchings of $G$ that are canonical in $G[A]$ and
non-canonical in $G[B]$ is at least $3|B|/4-7.5$.

Assume now that Case B applies for $i=1$ and $j=2$. The number of
matchings of $G$ that are canonical in $G[A]$ and non-canonical in
$G[B]$ is at least
\begin{equation}
\overline{m}^B_{1234} +\overline{m}^B_{13} +\overline{m}^B_{14}
+\overline{m}^B_{23} +\overline{m}^B_{24} +2\overline{m}^B_{\varnothing}.
\label{eq-4conn-1b}
\end{equation}
The number of perfect matching of $G^B_{13}$ is equal to the following
quantity which must be at least $|B|/2$ as argued before:
\begin{equation}
m^B_{1234}+m^B_{13}+m^B_{24}+m^B_{\varnothing}\ge |B|/2.
\label{eq-4conn-2b}
\end{equation}
Similarly, we bound the number of perfect matchings of $G^B_{14}$:
\begin{equation}
m^B_{1234}+m^B_{14}+m^B_{23}+m^B_{\varnothing}\ge |B|/2.
\label{eq-4conn-3b}
\end{equation}
Finally, we estimate the number of perfect matchings of $G^B_{(12)}$:
\begin{equation}
m^B_{13}+m^B_{14}+m^B_{23}+m^B_{24}+m^B_{\varnothing}\ge |B|/2+1.
\label{eq-4conn-4b}
\end{equation}
Summing (\ref{eq-4conn-2b}), (\ref{eq-4conn-3b}) and
(\ref{eq-4conn-4b}) and subtracting the maximum possible number of
canonical matchings, we obtain
\begin{equation}
2\overline{m}^B_{1234} +2\overline{m}^B_{13} +2\overline{m}^B_{14}
+2\overline{m}^B_{23} +2\overline{m}^B_{24}
+3\overline{m}^B_{\varnothing}\ge 3|B|/2-15.
\label{eq-4conn-5b}
\end{equation}
Comparing (\ref{eq-4conn-1b}) and (\ref{eq-4conn-5b}), we see that the
number of perfect matchings of $G$ that are canonical in $G[A]$ and
non-canonical in $G[B]$ is at least $3|B|/4-7.5$.

A completely symmetric argument yields that the number of perfect
matchings of $G$ that are non-canonical in $G[A]$ and canonical in
$G[B]$ is at least $3|A|/4-7.5$. We now consider matchings of $G$ that
are canonical when restricted to both $G[A]$ and $G[B]$. If Case A
applies to both $G[A]$ and $G[B]$, there are at least $6+2\cdot 2=10$
such perfect matchings of $G$. If Case A only applies to one of these
two subgraphs, there are at least $4+2\cdot 2=8$ such perfect
matchings. Finally, if Case B applies to both $G[A]$ and $G[B]$, there
are at least $2+2\cdot 2=6$ such perfect matchings. In total, the
number of perfect matchings of $G$ is at least
$3|A|/4-7.5+3|B|/4-7.5+6=3n/4-9$.
\end{proof}

\section{Cyclically $3$-edge-connected graphs}
\label{sect3}

A {\em klee-graph} is inductively defined as being either $K_4$, or
the graph obtained from a klee-graph by replacing a vertex by a
triangle. Every klee-graph is a cubic planar brick. Moreover, if $G$
is a graph with an edge-cut $E(A,B)$ such that both $G/A$ and $G/B$
are klee-graphs, then $G$ is also a klee-graph.

Recall that every edge of a cubic bridgeless graph is contained in at
least one perfect matching. We now prove that if an edge of a
$3$-edge-connected cubic graph is contained in only one perfect
matching, then the graph is a klee-graph.

\begin{lemma}
\label{lm-double}
A $3$-edge-connected cubic graph $G$ that is not a klee-graph is
matching double-covered.
\end{lemma}

\begin{proof}
The proof proceeds by induction on the order of $G$. If $G$ has no
cyclic $3$-edge-cuts, then it is matching double-covered by
Proposition~\ref{prop-double} (as $G$ is not a klee-graph, it is
different from $K_4$). Otherwise, let $E(A,B)$ be a cyclic
$3$-edge-cut of $G$. Since $G$ is not a klee-graph, at least one of
the graphs $G/A$ and $G/B$, say $G/A$, is not a klee-graph. By
induction, $G/A$ is matching double-covered. Since $G/B$ is cubic and
bridgeless, it is matching covered. Hence, every perfect matching of
$G/A$ extends to $G$, and so every edge with at least one end-vertex
in $B$ is contained in at least two perfect matchings of $G$.

If $e$ is an edge with both end-vertices in $A$, then there exists a
perfect matching of $G/B$ containing $e$. Since $G/A$ is matching
double-covered, this matching extends in two different ways to a
matching of $G$. Hence, $G$ is matching double-covered.
\end{proof}

\paragraph{}In this section, our general strategy to prove that a
cyclically 3-edge-connected cubic graph has many matchings is to split
the graph along a 3-edge-cut and then use an inductive argument. If
the smaller graphs are not klee-graphs, every edge of such graphs is
in at least two perfect matchings and those can be combined to form
many different matchings in the original graph.

\begin{lemma}
\label{lm-nonkubic}
Every $n$-vertex $3$-edge-connected cubic graph $G$ with a
$3$-edge-cut $E(A,B)$ such that neither $G/A$ nor $G/B$ is a klee-graph, has
at least $3n/4-6$ perfect matchings.
\end{lemma}

\begin{proof}
Let $E(A,B)=\{e_1,e_2,e_3\}$, and let $m^A_i$ (resp. $m^B_i$) be the
number of perfect matchings of $G/A$ (resp. $G/B$) containing the edge
$e_i$. By Lemma~\ref{lm-double}, each of $m^A_i$ and $m^B_i$ is at least
two. By Theorem~\ref{thm-previous}, unless $G/A$ is the exceptional
graph from Figure~\ref{fig-exceptional},
$$m^A_1+m^A_2+m^A_3\ge|B|/2+3/2\quad\mbox{and} \quad
m^B_1+m^B_2+m^B_3\ge|A|/2+1/2\;\mbox{.}$$ Since any perfect matching
of $G/A$ containing $e_i$ combines with a perfect matching of $G/B$
containing $e_i$ to form a perfect matchings of $G$ containing $e_i$,
the number of perfect matchings of $G$ is at least
$$\sum_{i=1}^3 m^A_im^B_i\ge 2\,(|B|/2-5/2)+2\,
(|A|/2-7/2)+2\cdot 2=|A|+|B|-8=n-8\;\mbox{.}$$ Since neither $G/A$ nor
$G/B$ is a klee-graph, and both $A$ and $B$ have odd size, $|A|\ge 5$ and
$|B|\ge 5$. Consequently, $n=|A|+|B|\ge 10$ and thus $G$ has at least
$n-8\ge 3n/4-5.5$ perfect matchings.

If $G/A$ is the exceptional graph, then $|B|=11$ and
$m^A_1=m^A_2=m^A_3=2$. The bound on the number of perfect matchings of
$G$ is now $$\sum_{i=1}^3 m^A_im^B_i\ge
2\,(|A|/2+1/2)=|A|+1=n-10\;\mbox{.}$$ Since $|B|=11$ and $|A|\ge 5$,
the number $n$ of vertices of $G$ is at least $16$, and so $G$ has at
least $n-10\ge 3n/4-6$ perfect matchings.
\end{proof}

We say that a $3$-edge-cut $E(A,B)$ of a cubic graph $G$ is
\emph{nice}, if $G/A$ is not a klee-graph and at least one of the following
holds:

  \smallskip
  \qite{(i)}$G/B$ is not a klee-graph;

  \smallskip
  \qite{(ii)}$|A|\ge 9$;

  \smallskip
  \qite{(iii)}$|A|\ge 5$ and $E(A,B)$ is not tight;

  \smallskip
  \qite{(iv)}$|A|=3$, and there are at least two perfect
      matchings of $G$ containing all the three edges of $E(A,B)$.\\

The next lemma shows that if we split the graph along a nice
3-edge-cut, the general induction will run smoothly.

\begin{lemma}
\label{lm-3cut}
Let $n$ be a positive integer, and assume that every
$3$-edge-connected cubic graph of order $n'<n$ has at least $3n'/4-9$
perfect matchings. If $G$ is an $n$-vertex $3$-edge-connected cubic
graph with a nice $3$-edge-cut $E(A,B)$, then $G$ also has at least
$3n/4-9$ perfect matchings.
\end{lemma}

\begin{proof}
By the assumption of the lemma, $G/A$ is not a klee-graph.  If $G/B$
is also not a klee-graph, the bound follows from
Lemma~\ref{lm-nonkubic}. We now focus on the remaining three cases and
assume that $G/B$ is a klee-graph. By Lemma~\ref{lm-double}, the graph
$G/A$ is matching double-covered.  Since $G/A$ has fewer vertices than
$G$, by the assumption of the lemma $G/A$ has at least $3|B|/4+3/4-9$
perfect matchings. Since $G/B$ is a klee-graph, we conclude that it is
not the exceptional graph from Figure~\ref{fig-exceptional}, and thus
it has at least $|A|/2+3/2$ perfect matchings.

Let $E(A,B)=\{e_1,e_2,e_3\}$, and let $m^A_i$ (resp. $m^B_i$) be the
number of perfect matchings of $G/A$ (resp. $G/B$) containing $e_i$,
$i=1,2,3$. The number of perfect matchings of $G$ containing exactly
one edge of the edge-cut $E(A,B)$ is at least
\begin{equation}
m^A_1\cdot m^B_1+ m^A_2\cdot m^B_2+ m^A_3\cdot m^B_3.
\label{eq-t1}
\end{equation}
As every $m^A_i$ is at least two and every $m^B_i$ is at least one,
the expression above is at least
\begin{equation}
(3|B|/4+3/4-13)\cdot 1+2\cdot (|A|/2-1/2)+2\cdot 1=%|A|+3|B|/4+3/4-12=
3n/4+|A|/4+3/4-12
\label{eq-t2}
\end{equation}
If $|A|\ge 9$, then $3n/4+|A|/4+3/4-12\ge 3n/4+12/4-12=3n/4-9$.  If
$|A|\ge 5$ and the edge-cut $E(A,B)$ is not tight, then there exists a
perfect matching not counted in the estimate (\ref{eq-t2}) and thus
the number of perfect matchings is at least $3n/4+|A|/4+3/4-11\ge
3n/4-9$. Finally, assume that $|A|=3$ and there are at least two
perfect matchings containing all the three edges of $E(A,B)$, i.e., at
least two matchings are not counted in (\ref{eq-t2}). Then the number of
perfect matchings of $G$ is at least $3n/4+|A|/4+3/4-10 > 3n/4-9$.
\end{proof}

Let $G$ and $H$ be two disjoint cubic graphs, $u$ a vertex of $G$
incident with three edges $e_1,e_2,e_3$, and $v$ a vertex of $H$
incident with three edges $f_1,f_2,f_3$. Consider the graph obtained
from the union of $G \backslash u$ and $H \backslash v$ by adding an
edge between the end-vertices of $e_i$ and $f_i$ ($1\le i \le3$)
distinct from $u$ and $v$. We say that this graph is obtained by
\emph{gluing} $G$ and $H$ through $u$ and $v$. Note that gluing a
graph $G$ and $K_4$ through a vertex $v$ of $G$ is the same as
replacing $v$ by a triangle.

In the next lemma, we characterize the graphs that do not contain nice
3-edge-cuts.

\begin{lemma}
\label{lm-special}
Let $G$ be a $3$-edge-connected cubic graph that is not cyclically
$4$-edge-connected and that has no nice $3$-edge-cut. If $G$ is
neither a klee-graph nor bipartite, then $G$ must be of one of the
following forms:

  \smallskip
  \qite{(1)}$G$ can be obtained from a cubic brace $H$ by gluing
      klee-graphs on $4$, $6$ or $8$ vertices through some of the
      vertices of one of the two color classes of $H$;

  \smallskip
  \qite{(2)}$G$ has no tight edge-cuts and can be obtained from a
      cyclically $4$-edge-connected cubic graph by replacing some of
      its vertices with triangles.

\end{lemma}

\begin{proof}
We assume that $G$ is neither a klee-graph nor a bipartite graph and
distinguish two cases depending whether $G$ has a tight edge-cut or
not.

If $G$ has a tight edge-cut, then its brick and brace decomposition is
non-trivial. Every non-trivial brick and brace decomposition of a
cubic bridgeless graph contains a brace (see~\cite{bib-previous}). If
the brick and brace decomposition of $G$ contains two or more braces,
then $G$ has a tight $3$-edge-cut $E(A,B)$ such that neither $G/A$ nor
$G/B$ is a brick (again, see~\cite{bib-previous}). In particular,
neither $G/A$ nor $G/B$ is a klee-graph, and so $E(A,B)$ is a nice
edge-cut, which violates the assumption of the lemma.

We conclude that the brick and brace decomposition of $G$ contains a
single brace $H$, and that for any tight edge-cut $E(A,B)$ of $G$,
exactly one of the graphs $G/A$ and $G/B$ is a brick. Observe that all
the bricks are glued through the vertices of the same color class of
$H$. To see this, assume that for two vertices $u$ and $v$ in
different color classes of $H$, and two bricks $H_1$ containing a
vertex $u'$ and $H_2$ containing a vertex $v'$, $G$ is obtained from
$H$ by gluing $H_1$ through $u$ and $u'$ and $H_2$ to $v$ and
$v'$. Let $u_1,u_2,u_3$ (resp. $v_1,v_2,v_3$) be the neighbors of $u$
(resp. $v$) in $H$, and let $u'_1,u'_2,u'_3$ (resp. $v'_1,v'_2,v'_3$)
be the neighbors of $u'$ (resp. $v'$) in $H_1$ (resp. $H_2$). By
definition, both $\{u_iu'_i, 1\le i \le 3\}$ and $\{v_iv'_i, 1\le i
\le 3\}$ are tight edge-cuts of $G$. Since $H_1$ and $H_2$ are bricks,
$H_1\backslash \{u'_1,u'_2\}$ and $H_2\backslash \{v'_1,v'_2\}$ both
have a perfect matching. Since $H$ is a brace, $H \backslash
\{u_1,u_2,v_1,v_2\}$ also has a perfect matching. These three
matchings combine to a perfect matching of $G$ containing all the edges
$u_iu'_i$ and $v_iv'_i$ for $1 \le i \le 3$ which contradicts the
fact that the two edge-cuts were tight.

As for every $3$-edge-cut $E(A,B)$, $G/A$ or $G/B$ is a klee-graph,
all bricks of $G$ are klee-graphs. Since $E(A,B)$ is not nice, the
``klee-graph'' side of the cut has at most $8$ vertices. Hence, all
bricks of $G$ are klee-graphs with $4$, $6$ or $8$ vertices, and $G$
is exactly of the first form described in the lemma.

It remains to consider the case that $G$ has no tight $3$-edge-cuts.
Consider a 3-edge-cut $E(A,B)$ of $G$. Since $G$ is not a klee-graph,
$G/A$ or $G/B$, say $G/A$, is not a klee-graph. Since $G$ has no nice
$3$-edge-cut, $|A|=3$ and so $G[A]$ is a triangle. Now observe that
every 3-edge-cut in $G/A$ corresponds to a 3-edge-cut in $G$, and
hence, separates a triangle. So we can keep contracting the original
triangles of $G$ to obtain a cyclically 4-edge-connected graph (no new
3-edge-cut, and hence no triangle, will be created during the
process). We have observed that $G$ can be obtained from a cyclically
$4$-edge-connected cubic graph by replacing some of its vertices by
triangles.
\end{proof}

Let $G$ be a $3$-edge-connected cubic graph that is not a klee-graph,
such that every cyclic $3$-edge-cut $E(A,B)$ of $G$ separates a
triangle (in other words $|A|=3$ or $|B|=3$). The {\em core} of $G$,
denoted by ${\cal C}(G)$, is the graph obtained by contracting every
triangle of $G$. Since all cyclic $3$-edge-cuts of $G$ separate
triangles, the graph $G$ can be obtained from its core by replacing
some of its vertices with triangles.

\begin{lemma}
\label{lm-special2}
Let $G$ be a $3$-edge-connected cubic graph different from $K_4$ with
no nice $3$-edge-cut. Assume $G$ was obtained from a cyclically
$4$-edge-connected cubic graph by replacing some of its vertices (at
least one) by triangles. In particular, $G$ is not a klee-graph. If ${\cal
C}(G)$ is not bipartite, then ${\cal C}(G)$ has a cyclic $4$-edge-cut,
and $G$ has no tight cyclic $3$-edge-cut.
\end{lemma}

\begin{proof}
Let $H={\cal C}(G)$ and let $v$ be any vertex of $H$. By the
assumption, $H$ is not bipartite. If the graph $H'$ obtained from $H$
by removing $v$ and its three neighbors has no perfect matching, then
there exists $S'\subseteq V(H')$ such that $H\setminus S'$ has at
least $|S'|+2$ odd components. Let $S$ be the set $S'$ enhanced with
the three neighbors of $v$. Clearly, $H\setminus S$ has at least
$|S|=|S'|+3$ odd components.  Since $H$ is cyclically
$4$-edge-connected, this implies that all the odd components of
$H\setminus S$ are isolated vertices and $H$ is bipartite which is
impossible. Hence, $H'$ has a perfect matching.

 Let $u$ be a vertex of $H$ that is replaced by a triangle $T$ in $G$
and let $U$ be the set containing $u$ and its three neighbors
$u_1,u_2, u_3$ in $H$. As proven in the previous paragraph,
$H\setminus U$ contains a perfect matching and the cut separating the
triangle $T$ is not tight. Hence, no cyclic $3$-edge-cut of $G$ is
tight.

We now show that $H$ has a cyclic 4-edge-cut. If $H\setminus U$
contains two perfect matchings, then $G$ has two perfect matchings
containing all the three edges of the cut separating $T$. Since $G$
has no nice $3$-edge-cut, this is impossible, so by Kotzig's theorem
the graph $H\setminus U$ has a bridge. Let $E(A,B)$ be the cut of
$H\setminus U$, that corresponds to this bridge.

Since $H$ is cyclically $4$-edge-connected, the set $\{u_1,u_2,u_3\}$
is a stable set. If $A$ is comprised of a single vertex, say
$A=\{v\}$, then $v$ has two common neighbors with $u$, say $u_1$ and
$u_2$. In particular, $H$ contains the cycle of length four $uu_1vu_2$
which is disjoint from $B$. If $B$ induces a forest it is easy to see
that $|B|=3$ and $B$ induces a path of length two, which together with
$u_3$ forms a cycle of length four. Otherwise, $B$ has a cycle. In
both cases, $H$ has a cyclic edge-cut of size four. Since the case
$|B|=1$ is symmetric, we can assume that both $A$ and $B$ contain at
least two vertices. Since $H$ is cyclically $4$-edge-connected, the
sizes of the cuts $E(A,B\cup U)$ and $E(A\cup U,B)$ are at least
four. Since the number of edges between $U$ and $A\cup B$ is six,
there are three edges joining $U$ and $A$ and three edges joining $U$
and $B$.

If $|A|\ge 3$, then $E(A,B\cup U)$ is a cyclic edge-cut of size
four. If $|A|=2$, then one of the two vertices of $A$ has two common
neighbors with $u$ and $H$ has a cycle of length four. Again, $H$ has
a cyclic edge-cut of size four.
\end{proof}

As mentioned in the introduction, Chudnovsky and
Seymour~\cite{bib-chudnovsky08+} proved that planar cubic bridgeless
graphs (and consequently, klee-graphs) have exponentially many
perfect matchings. However, their bound is not too good for graphs
with small number of vertices. In the next lemma, we use the inductive
structure of klee-graphs to provide a better lower bound on their
number of perfect matchings.

\def\minord{10}
\def\subst{6}

\begin{lemma}
\label{lm-kubic}
Every $n$-vertex klee-graph has at least $3n/4-\subst$ perfect
matchings.
\end{lemma}

\begin{proof}
If $n \le 8$, then there is nothing to prove. Hence,
we can focus on klee-graphs of order at least ten.

Let $G$ be a klee-graph and $v$ a vertex of $G$ with neighbors $v_1$,
$v_2$ and $v_3$. The \emph{type} of $v$ is the 4-tuple
$(\omega;\mu_1,\mu_2,\mu_3)$ such that the graph
$G\setminus\{v,v_1,v_2,v_3\}$ contains $\omega$ perfect matchings and
the graph $G\setminus\{v,v_i\}$ contains $\mu_i$ perfect matchings for
$1 \le i \le 3$. Observe that there are exactly three non-isomorphic
klee-graphs of order ten; these graphs are depicted in
Figure~\ref{fig:kubic}(a)--(c), where the label of each edge
represents the number of perfect matchings containing that edge and
the label of a vertex $v$ is the number of perfect matchings in the
graph obtained by removing $v$ and its three neighbors. In particular,
the type of a vertex $v$ is formed by its label and the labels
of the three incident edges.

\begin{figure}[htbp]
\begin{center}
\hspace{0.05cm}
\subfigure[\label{fig:k1}]{\includegraphics[scale=0.75]{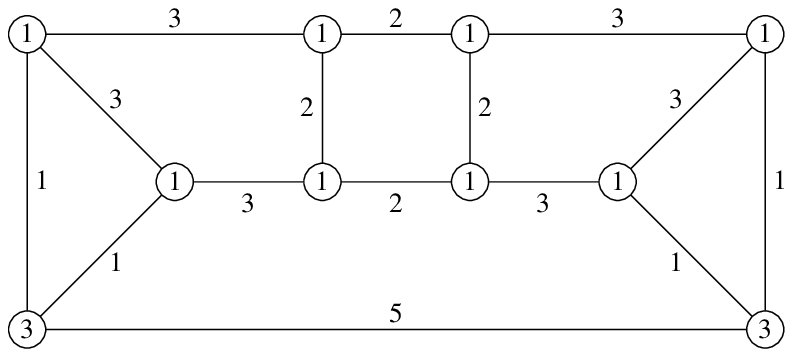}}
\hspace{0.6cm}
\subfigure[\label{fig:k2}]{\includegraphics[scale=0.75]{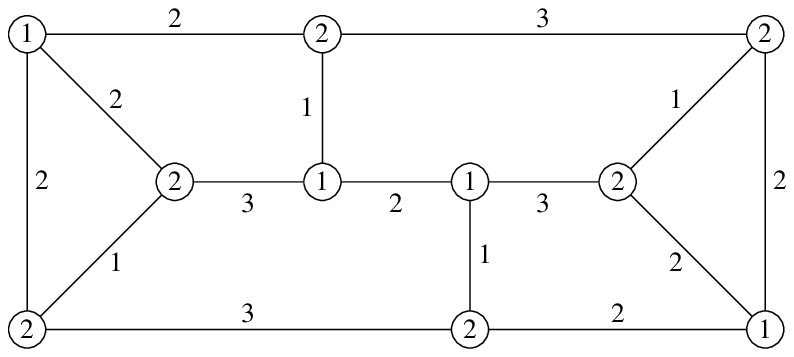}}\vskip 5mm
\subfigure[\label{fig:k3}]{\includegraphics[scale=0.75]{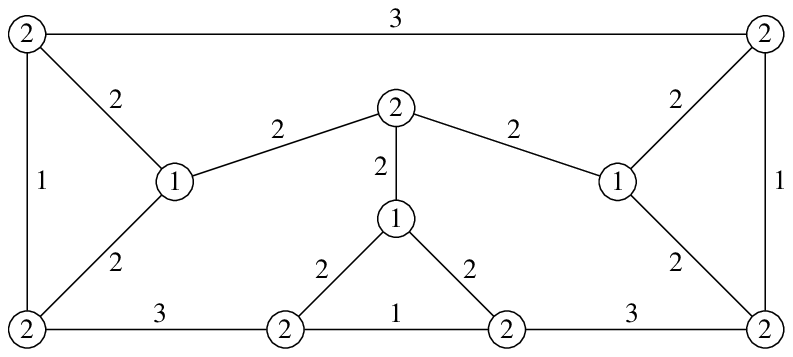}}
\hspace{0.6cm}
\subfigure[\label{fig:k4}]{\includegraphics[scale=0.75]{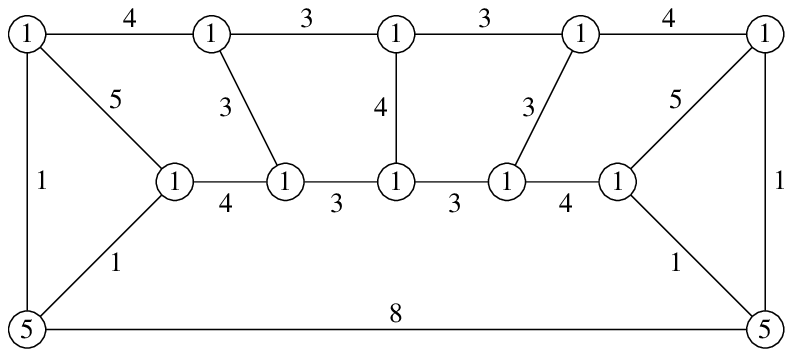}}
\caption{(a)--(c) The three non-isomorphic klee-graphs of order ten. (d)
         The only $12$-vertex klee-graph that cannot be obtained by
         replacing a vertex by a triangle in (b) or (c). %% The label
%%          of each edge represents the number of perfect matchings
%%          containing that edge and the label of a vertex $v$ is the
%%          number of perfect matchings in the graph obtained by removing
%%          $v$ and its three neighbors. In particular, the type of a
%%          vertex $v$ is formed by its label and the labels of the three
%%          incident edges.
	 \label{fig:kubic}}
\end{center}
\end{figure}

%% \begin{figure}
%% \begin{center}
%% \epsfbox{k1}\vskip 5mm
%% \epsfbox{k2}\vskip 5mm
%% \epsfbox{k3}
%% \end{center}
%% \caption{The three non-isomorphic klee-graphs of order ten. The label of
%%          each edge represents the number of perfect matchings containing
%% 	 that edge and the label of a vertex $v$ is the number of perfect
%% 	 matchings in the graph obtained by removing $v$ and its three
%% 	 neighbors. In particular, the type of a vertex $v$ is formed
%% 	 by its label and the labels of the three incident edges.}
%% \label{fig-kubic10}
%% \end{figure}

Let $v$ be a vertex of type $(\omega;\mu_1,\mu_2,\mu_3)$ in the
klee-graph $G$. The vertex $v$ is said to be an \emph{$A$-vertex} if
$\omega=1$ and $\mu_i=1$ for a single index $i\in\{1,2,3\}$; $v$ is a
\emph{$B$-vertex} if $\omega=1$ and $\mu_i>1$ for every
$i\in\{1,2,3\}$ and $v$ is a \emph{$C$-vertex} if $\omega>1$ and
$\mu_i=1$ for exactly two indices $i\in\{1,2,3\}$. A vertex is
\emph{dangerous} if at least three of the values $\omega$, $\mu_1$,
$\mu_2$ and $\mu_3$ are equal to one. A vertex $v$ is \emph{good} if
it is neither a $A$-, $B$-, $C$-vertex nor a dangerous vertex. In the
following, $G \! \bigtriangleup \! v$ denotes the graph obtained from
$G$ by replacing $v$ with a triangle. The number of perfect matchings
in $G$ is denoted by $m(G)$.

Let $G$ be a klee-graph and $v$ a vertex of $G$ of type
$(\omega;\mu_1,\mu_2,\mu_3)$. As illustrated in
Figure~\ref{fig:exkubic}, the types of the three new vertices in $G \!
\bigtriangleup \! v$ are
$$(\mu_1;\mu_1+\omega,\mu_2,\mu_3)\mbox{, }
(\mu_2;\mu_1,\mu_2+\omega,\mu_3)\mbox{, and }
(\mu_3;\mu_1,\mu_2,\mu_3+\omega)\;.$$ In particular, $m(G \!
\bigtriangleup \!  v)=m(G)+ \omega$.  Finally, consider a vertex
$v'\not=v$ and observe that if the type of $v'$ in $G$ is
$(\omega';\mu'_1,\mu'_2,\mu'_3)$ and its type in $G \! \bigtriangleup
\! v$ is $(\omega'';\mu''_1,\mu''_2,\mu''_3)$, then
$\omega''\ge\omega'$ and $\mu''_i\ge\mu'_i$ for every
$i\in\{1,2,3\}$. Hence, if $v'$ is an $A$-vertex in $G$, it is an
$A$-vertex, a $B$-vertex or a good vertex in $G\! \bigtriangleup \!
v$. If $v'$ is a $B$-vertex in $G$, it is a $B$-vertex or a good
vertex in $G\! \bigtriangleup \! v$. If $v'$ is a $C$-vertex in $G$,
then it is a $C$-vertex or a good vertex in $G\!  \bigtriangleup \!
v$. Finally, if $v'$ is a good vertex in $G$, it remains
good in $G\!  \bigtriangleup \! v$. This implies that a vertex is
dangerous in $G \!  \bigtriangleup \!  v$ only if it was dangerous in
$G$. Since no graph in Figure~\ref{fig:kubic}(a)--(c) contains a
dangerous vertex, no klee-graph of order at least 12 contains a
dangerous vertex.

\begin{figure}[htbp]
\begin{center}
\includegraphics[scale=0.9]{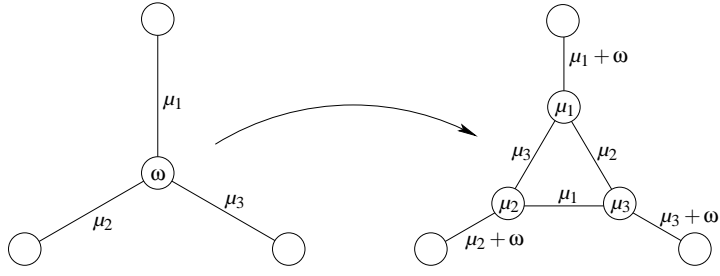}
\caption{The types of the three new vertices in $G \! \bigtriangleup \!
v$. \label{fig:exkubic}}
\end{center}
\end{figure}

For any klee-graph $G$ with $\alpha$ $A$-vertices and $\beta$
$B$-vertices, let $M(G)=m(G)-\alpha-\beta/2$. The core of our proof is
the following claim proven by induction on $n$.\\

\noindent {\bf Claim.} {\em For any $n$-vertex klee-graph $G$,
$n\ge\minord$, distinct from the one in Figure~\ref{fig:k1}, it holds
$M(G) \, \ge \, 3n/4-\subst$.}\\

\noindent If $n=\minord$, then $G$ is one of the graphs depicted in
Figures~\ref{fig:k2} and~\ref{fig:k3}, and
$$M(G)=\left\{\begin{array}{ccc}
       %7-4-3\cdot 2/2-4/2 & = & -2 \\
       6-2-2/2 & = & 3 \\
       6-0-3/2 & = & 4.5
       \end{array}
       \right\}\ge3\cdot 10/4-\subst\;\mbox{.}$$ The only $12$-vertex
klee-graph that cannot be obtained by replacing a vertex with a
triangle in one of the graphs depicted in Figures~\ref{fig:k2}
and~\ref{fig:k3} is the graph in Figure~\ref{fig:k4}. For this graph,
we have
$$M(G)=10-4-6/2=3\ge 3\cdot 12/4-\subst\;\mbox{.}$$ All other
$n$-vertex klee-graphs $G$ with $n\ge 12$ can be obtained by
replacing a vertex $v$ by a triangle $w_1w_2w_3$ in a klee-graph $G'$
that satisfies the assumptions of the claim.  Clearly, the number $n'$
of vertices of $G'$ is $n-2$. By the induction, we assume that
$M(G')\ge 3n'/4-\subst$.

%% \begin{figure}
%% \begin{center}
%% \epsfbox{k4}
%% \end{center}
%% \caption{The only $12$-vertex klee-graph that cannot be obtained
%%          by replacing a vertex with a triangle from the second or
%% 	 the third graph depicted in Figure~\ref{fig-kubic10}.}
%% \label{fig-kubic12}
%% \end{figure}

We now distinguish four cases based on the type of $v$; note that $v$
cannot be dangerous as argued earlier. Observe that if an $A$- or
$B$-vertex becomes good, or if an $A$-vertex becomes a $B$-vertex,
then $-\alpha-\beta/2$ increases. So we can assume without loss of
generality that every $A$-vertex and $B$-vertex distinct from $v$
remains an $A$-vertex and $B$-vertex, respectively.
\begin{description}
\item[$\bullet$] \emph{$v$ is an $A$-vertex:} Since $v$ is
     an $A$-vertex, $m(G)=m(G')+1$.  One of the vertices $w_1$, $w_2$
     and $w_3$ is a $B$-vertex, and the other two vertices are
     good. Hence, $\alpha$ decreases by $1$ and $\beta$ increases by
     $1$, and so $-\alpha-\beta/2$ increases by $1/2$. We conclude
     that
     $$M(G) \, = \, M(G')+1+1/2 \, \ge \, 3n'/4-\subst+3/2 \, = \,
     3n/4-\subst\;\mbox{.}$$
\item[$\bullet$] \emph{$v$ is a $B$-vertex:} Since $v$ is a
     $B$-vertex, it holds that $m(G)= m(G')+1$. All the vertices
     $w_1$, $w_2$ and $w_3$ are good, so $\beta$ decreases by one and
     $-\alpha-\beta/2$ increases by $1/2$. We conclude that
     $$M(G) \,=  \,M(G')+1+1/2 \,\ge \, 3n'/4-\subst+3/2 \,= \,3n/4-\subst\;\mbox{.}$$
\item[$\bullet$] \emph{$v$ is a $C$-vertex:} It is easy to
     see that in any klee-graph of order at least 12, any $C$-vertex
     has type $(\omega,\mu,1,1)$, where both $\omega$ and $\mu$ are at
     least five. Hence it holds that $m(G)\ge m(G')+5$. Two vertices
     among $w_1$, $w_2$ and $w_3$ are $A$-vertices and the last one is
     a $C$-vertex. Hence, $-\alpha-\beta/2$ decreases by two. We again
     conclude that
     $$M(G) \,\ge \, M(G')+5-2 \, \ge \, 3n'/4-\subst+3 \,\ge \,3n/4-\subst\;\mbox{.}$$
\item[$\bullet$] \emph{$v$ is good:} At most one of the
     vertices $w_1$, $w_2$ and $w_3$ is a $B$-vertex and the remaining
     vertices are good. Hence, $-\alpha-\beta/2$ decreases by
     at most $1/2$. Since $m(G)\ge m(G')+2$, it holds that
     $$M(G) \, \ge \, M(G')+2-1/2 \, \ge \, 3n'/4-\subst+3/2 \, =
     \,3n/4-\subst\;\mbox{.}$$
\end{description}
This finishes the proof of the claim.

We have shown that $M(G)\ge 3n/4-\subst$ for every $n$-vertex
klee-graph $G$ with $n\ge\minord$ distinct from the graph in
Figure~\ref{fig:k1} which has $7 \ge 3\cdot10/4-6$ perfect
matchings. In particular, the number of perfect matchings of any
$n$-vertex klee-graph is at least $3n/4-\subst$.
\end{proof}

As mentioned in the introduction, cubic bridgeless bipartite graphs
are known to have an exponential number of perfect matchings. We can
derive the following more modest result, which will be sufficient for
our purpose.

\begin{lemma}
\label{lm-bipartite}
Every $n$-vertex cubic bipartite graph has at least $3n/2-9$ perfect
matchings.
\end{lemma}

\begin{table}
\begin{center}
\begin{tabular}{|c|c|c|c|c|c|c|}
\hline
  $n=2k$ &  6 &  8 & 10 & 12 & 14 & 16 \\
\hline
  $g(k)$ &  4 &  6 &  8 & 11 & 15 & 20 \\
\hline
  $f(k)$ &  6 &  9 & 12 & 17 & 23 & 30 \\
\hline
$3n/2-9$ &  0 &  3 &  6 &  9 & 12 & 15 \\
\hline
\end{tabular}
\end{center}
\caption{The minimum number $f(k)$ of distinct perfect matchings of a
         cubic bipartite with $2k$ vertices and the claimed bound
         $3n/2-9$.}
\label{table-bip}
\end{table}

\begin{proof}
Let $g(3)=4$, and set $g(k)=\left\lceil 4g(k-1)/3\right\rceil$ for any
$k\ge 4$. Also, let $f(k)=\left\lceil 3g(k)/2\right\rceil$.  It can be
shown that every cubic bridgeless graph with $2k$ vertices has at
least $f(k)$ perfect matchings, see~\cite{bib-previous,bib-lovasz86+}.
The values of $f(k)$ for small $k$ can be found in
Table~\ref{table-bip}. If $n\le 12$, the statement of the lemma holds
by inspecting the values of $f(k)$. For $k=7$, $g(k)\ge 2k$. Using the
definition of $g(k)$, an easy argument by induction on $k$ shows that
$g(k)\ge 2k$ for all $k\ge 7$.  Hence, $f(k)\ge 3g(k)/2\ge 3k=3n/2$
and the statement of the lemma follows.
\end{proof}

We are know ready to prove the main result of this section.

\begin{theorem}
\label{thm-3conn}
Every $n$-vertex $3$-edge-connected cubic graph has at least $3n/4-9$
perfect matchings.
\end{theorem}

\begin{proof}
The proof proceeds by induction on the order $n$ of $G$. If $n\le 12$,
then there is nothing to prove since the bound claimed in the theorem
is negative. Fix $n \ge 14$, and assume that we have proven the
statement of the theorem for all $n'<n$. If $G$ is cyclically
$4$-edge-connected, then $G$ has at least $3n/4-9$ perfect matchings
by Theorem~\ref{thm-4conn}. If $G$ has a nice cyclic $3$-edge-cut,
then $G$ has at least $3n/4-9$ perfect matchings using
Lemma~\ref{lm-3cut}. If $G$ is a klee-graph or a bipartite graph,
Lemmas~\ref{lm-kubic} and~\ref{lm-bipartite} yield the desired lower
bound on the number of perfect matchings of $G$.  Otherwise, $G$ is of
one of the two forms given in Lemma~\ref{lm-special}. We deal with
each of these cases separately:
\begin{itemize}
\item \emph{$G$ can be obtained from a cubic brace $H$ by gluing
      klee-graphs on $4$, $6$ or $8$ vertices through some of the
      vertices of one of the two color classes of $G$:} Let $N$ be the
      order of $H$. The number of perfect matchings of $H$ is at least
      $3N/2-9$ by Lemma~\ref{lm-bipartite} and $H$ is matching
      double-covered by Lemma~\ref{lm-double}. Let $N_k$ be the number
      of vertices of $H$ through which a klee-graph of order
      $k\in\{4,6,8\}$ is glued.  Observe that $$N_4+N_6+N_8\le N/2
      \mbox{ and } n=N+2N_4+4N_6+6N_8\;.$$

      Let us estimate the number of perfect matchings of $G$ in more
      detail. We count in how many ways perfect matchings of $H$
      extend to the glued klee-graphs. There is a unique extension of
      each perfect matching of $H$ to a glued klee-graph of order
      $4$. Since the edges incident with every vertex of a klee-graph
      of order six are contained in $1$, $1$ and $2$ perfect matchings
      respectively and $H$ is matching double-covered, at least two
      perfect matchings extend to a glued klee-graph of order six in
      two different ways. Hence, any such gluing increases the
      number of perfect matchings by at least two. Similarly, the
      edges incident with every vertex of a klee-graph of order eight
      are contained in $1$, $1$ and $3$ or $1$, $2$ and $2$ perfect
      matchings which implies that at least two matchings of $H$
      extend to a glued klee-graph of order eight in three different
      ways or at least four matchings of $H$ extend in two different
      ways.  In both cases, the number of perfect matchings is
      increased by four.

      Using Lemma~\ref{lm-bipartite}, we conclude that the
      number of perfect matchings of $G$ is at least
      $$\begin{array}{rcl} \frac32N-9+2N_6+ 4N_8\! & \!\ge\!  & \!
         \frac34N+3\,(N_4\!+\!N_6\!+\!N_8)/2+2N_6+4N_8-9 \\ & \!\ge\!
         & \! 3n/4-9\;,
       \end{array}$$ as desired.
\item \emph{$G$ has no tight edge-cuts and it can be obtained from a
      cyclically $4$-edge-connected cubic graph by replacing some of
      its vertices with triangles:} If $H={\cal{C}}(G)$ has a cyclic
      $4$-edge-cut, Theorem~\ref{thm-4conn} yields the desired bound.
      If $H$ has no cyclic $4$-edge-cut, then $H$ is a bipartite
      cyclically $5$-edge-connected cubic graph by
      Lemma~\ref{lm-special2}. By Proposition~\ref{prop-tight}, $H$ is
      a brace. In particular, it is possible to remove two vertices
      from each of the two colors classes of $H$ and the graph still
      has a perfect matching.

      Let $N$ be the number of vertices of $H$ and $N_i$, $i=1,2$, be
      the number of vertices of each of the two color classes of $H$
      that are replaced by triangles in $G$. Observe that
      $n=N+2N_1+2N_2$, $N_1\le N/2$ and $N_2\le N/2$. We can assume
      without loss of generality that $1 \le N_1\le N_2$, since
      otherwise this would bring us to the previous case (replacing
      a vertex $v$ by a triangle is the same as gluing a $K_4$ through
      $v$).

By Lemma~\ref{lm-bipartite}, $H$ has at least $3N/2-9$ perfect
      matchings and each of these matchings corresponds to a perfect
      matching of $G$ which contains only one edge of each 3-edge-cut
      separating a triangle.  Now, take two vertices $u,v$ in
      different color classes of $H$, such that $u$ and $v$ are
      replaced by two triangles $T_u$ and $T_v$ in $G$. Let $H'$ be
      the graph obtained from $H$ by removing two neighbors of $u$ and
      two neighbors of $v$. Since $H$ is a brace, $H'$ has a perfect
      matching. This perfect matching corresponds to a perfect
      matching of $G$ containing the three edges leaving $T_u$, the
      three edges leaving $T_v$, and only one edge of each 3-edge-cut
      separating a different triangle. Hence, $G$ contains at least
      $3N/2 -9 +N_1N_2$ perfect matchings. 

Since $n=N+2N_1+2N_2$,
      proving that $G$ has at least $3n/4-9$ perfect matchings is
      equivalent to proving that $N_1+N_2 \le \frac{N}2 + \frac23 N_1
      N_2$. If $N_1=1$ then $$N_1+N_2 = N_2/3 + 1 + \tfrac23 N_1N_2 \le
      N/2 + \tfrac23 N_1N_2$$ since $N \ge \lceil n/3 \rceil \ge 5$. On
      the other hand, if $N_1 \ge 2$ then $$N_1+N_2 \le N/2 +
      (N_1+N_2)/2 \le N/2 + N_1N_2/2\;.$$
\end{itemize}
This finishes the proof of Theorem~\ref{thm-3conn}.
\end{proof}

\section{Bridgeless graphs}
\label{sect2}

In this section, we prove our main result on the number of perfect
matchings of cubic bridgeless graphs. Before we do so, we need an
auxiliary lemma:

\begin{lemma}
\label{lm-avoid}
Let $G$ be a cubic bridgeless graph with a $2$-edge-cut. For every
edge $e$ of $G$, there are at least three perfect matchings avoiding
$e$.
\end{lemma}

\begin{proof}
Let $E(A,B)$ be an edge-cut of $G$ of size two and let $G^A$
and $G^B$ be the cubic bridgeless graphs obtained from $G[A]$
and $G[B]$ by joining the two vertices of degree two with an
edge. The added edges are denoted by $e^A$ and $e^B$. If $e
\in E(A,B)$, then $G$ has at least four perfect matchings avoiding $e$
as any of at least two perfect matchings of $G^A$ avoiding $e^A$
combines with any of at least two perfect matchings of $G^B$ avoiding
$e^B$ to a perfect matching of $G$ avoiding $e$.

We now assume that $e\not\in E(A,B)$. By symmetry, let $e$ be in
$G[A]$. Recall that in a cubic bridgeless graph, it is possible to
find a perfect matching avoiding any two given edges. Thus, the graph
$G^A$ contains at least two perfect matchings avoiding $e$ and at
least one such matching also avoids $e^A$. Any perfect matching of
$G^A$ avoiding both $e$ and $e^A$ can be extended to $B$ in two
different ways and any perfect matching of $G^A$ avoiding $e$ and
containing $e^A$ can be extended to $B$ in at least one
way. Altogether, $G$ contains at least three perfect matchings
avoiding $e$ as desired.
\end{proof}

We are now ready to prove the main result:

\begin{theorem}
\label{thm-main}
Every cubic bridgeless graph $G$ with $n$ vertices has at least
$3n/4-10$ perfect matchings.
\end{theorem}

\begin{proof}
The proof proceeds by induction on the number of vertices of $G$.  If
$G$ is $3$-edge-connected, the bound follows from
Theorem~\ref{thm-3conn}. Otherwise, take a $2$-edge-cut $E(A,B)$ of
$G$ such that $A$ is minimal with respect to inclusion.  Let $G^A$ and
$G^B$ be the cubic bridgeless graph obtained from $G[A]$ and $G[B]$ by
adding edges $e^A$ and $e^B$ between the two vertices of degree
two. Clearly, $G^A$ is 3-edge-connected and contains at least
$3|A|/4-9$ perfect matchings by Theorem~\ref{thm-3conn}. Also note
that $G^A$ contains at least two perfect matchings avoiding $e^A$, and
similarly $G^B$ contains at least two perfect matchings avoiding
$e^B$.

Suppose first that the edge $e^A$ is contained in two perfect
matchings. Fix two perfect matchings of $G^A$ containing $e^A$ and two
perfect matchings avoiding $e^A$.  Each of $|B|/2$ perfect matchings
of $G^B$ can be extended to $G[A]$ in at least two different ways
using the fixed matchings (note that, by Theorem~\ref{thm-previous},
if $|B|\not=12$, $G^B$ has at least $|B|/2+1$ perfect matchings and if
$|B|=2$, $G^B$ has $|B|/2+2=3$ perfect matchings). On the other hand,
every of at least $3|A|/4-9-4=3|A|/4-13$ perfect matchings of $G^A$
distinct from the fixed ones can be extended to $G[B]$.  Hence, unless
$|B|=2$ or $|B|=12$ the number of perfect matchings of $G$ is at least
$$3|A|/4-13+2\cdot(|B|/2+1)\, =\, 3n/4+|B|/4-11 \, \ge \, 3n/4-10\;.$$
If $|B|=2$, the number of perfect matchings of $G$ is at least
$$3|A|/4-13+2\cdot(|B|/2+2)\ge 3n/4-9\;,$$ and if $|B|=12$, the number
of perfect matchings of $G$ is at least $$3|A|/4-13+2\cdot |B|/2\, =\,
3n/4+|B|/4-13 \,= \, 3n/4-10\;.$$

Suppose now that $G^A$ has a single matching containing the edge
$e^A$. We distinguish two cases regarding whether $G^B$ is
$3$-edge-connected. If $G^B$ is $3$-edge-connected and $e^B$ is
contained in at least two perfect matchings, then we apply the same
arguments as in the previous paragraph and the result follows. Hence,
we can assume that $e^B$ is contained in a single perfect matching of
$G^B$. Consequently, by Theorem~\ref{thm-previous} there are at least
$|A|/2-1$ perfect matchings of $G^A$ avoiding $e^A$ and at least
$|B|/2-1$ perfect matchings of $G^B$ avoiding $e^B$. Fix two matchings
of $G^A$ that avoid $e^A$ and two matchings of $G^B$ that avoid $e^B$,
and call these four matchings \emph{canonical}. Every non-canonical
matching of $G^A$ avoiding $e^A$ combines with a canonical matching of
$G^B$ avoiding $e^B$, and vice-versa. Hence, the number of perfect
matchings of $G$ is at least $$2(|A|/2-3)+2(|B|/2-3)+2\cdot 2\, =\,
n-8 \, \ge \, 3n/4-9\;.$$

The only remaining case is when $G^B$ is not $3$-edge-connected and
the edge $e^A$ is contained in a single matching of $G^A$. By
Lemma~\ref{lm-avoid}, $G^B$ has at least three matchings avoiding
$e^B$. Fix one matching of $G^A$ containing $e^A$, one matching of
$G^A$ avoiding $e^A$ and three matchings of $G^B$ avoiding
$e^B$. Again, we call these five perfect matchings \emph{canonical}.
By induction, $G^B$ has at least $3|B|/4-10$ perfect matchings, each
of which can be combined with a canonical perfect matching of $G^A$ to
form a perfect matching of $G$. Since $e^A$ is contained in a single
matching of $G^A$, there exist at least $|A|/2-2$ matchings of $G^A$
(distinct from the canonical ones) avoiding $e^A$. Each of them can be
combined with one of the three canonical matchings of $G^B$ to form a
perfect matching of $G$. Note that $|A|/2-2\ge |A|/4$ if $|A|\ge 8$.
If $|A|\in\{4,6\}$, then by Theorem~\ref{thm-previous}, $G^A$ has at
least $|A|/2-1$ matchings distinct from the two canonical ones, and
again $|A|/2-1\ge |A|/4$. Finally, if $|A|=2$, then $G^A$ has
$|A|/2=1$ perfect matching distinct from the two canonical ones. In
all cases, $G^A$ has at least $|A|/4$ perfect matchings distinct from
the two canonical matchings of $G^A$. We conclude that the number of
perfect matchings of $G$ is at least
$$3\cdot|A|/4+3|B|/4-10\, =\, 3n/4-10\;.$$
This finishes the proof of the theorem.
\end{proof}

\end{document}